\documentclass[11pt]{article}

\usepackage{amsmath,amsthm,amsfonts}

\usepackage{amssymb}

\usepackage{epsfig}

\usepackage{rotating}

\usepackage{subfigure}

\newcommand{\p}{\partial}

\begin{document}

\title{Model reduction and mesh refinement}
\author{Panos Stinis \\ 
Department of Mathematics \\
University of Minnesota \\
    Minneapolis, MN 55455} 
%\address{Department of Mathematics, 
%    University of California 
%and Lawrence Berkeley National Laboratory, 
%    Berkeley, CA 94720}
%\email{stinis@math.lbl.gov}    
\date {}
%\thanks{This work was supported in part by the Director, Office of Science, Advanced Scientific Computing Research, U.S. Department of Energy under Contract No.
%DE-AC02-05CH11231.}

\maketitle

\begin{abstract}
In \cite{S07} we presented a novel algorithm for mesh refinement which utilizes a reduced model. In particular, the reduced model is used to monitor the transfer of activity (e.g. mass, energy) from larger to smaller scales. When the activity transfer rate exceeds a prescribed tolerance, the algorithm refines the mesh accordingly. The algorithm was applied to the inviscid Burgers and focusing Schr\"odinger equations to detect singularities. We offer here a simple proof of why the use of a reduced model for mesh refinement is a good idea. We do so by showing that by controlling the transfer of activity rate one controls the error caused inevitably by a finite resolution. 
\end{abstract}

%\maketitle

\section*{Introduction}

The subjects of mesh refinement and model reduction have attracted significant attention in the applied mathematics research community in the last couple of decades (see e.g. \cite{almgren,berger1,berger2,budd,ceniceros,CHK00,CS05,givon,landman}). The reason is that despite the ever increasing computational power there still exists a wealth problems that are not amenable to direct numerical simulation. The purpose of the current work is to establish a correspondence between these two seemingly disparate subjects. On the one hand, mesh refinement is concerned with the judicious allocation of computational power to allow focusing on the important characteristics of a solution. On the other hand, model reduction is concerned with the best way to approximate the effect of the unresolved scales of a simulation to the resolved ones. The element that connects these two subjects is the observation set forth in \cite{S07} that successful mesh refinement relies on the accurate monitoring of transfer of activity to the unresolved scales which is exactly the objective a good reduced model. Thus, if one possesses a good reduced model, it can be utilized to decide when and where to perform mesh refinement.

In \cite{S07} in order to determine the need to perform mesh refinement a reduced model was constructed for a subset of the resolved scales. Then, the rate of change of the $L_2$ norm predicted by the reduced model for this subset of the resolved scales was monitored (note that depending on the physical context the $L_2$ norm may represent, for example, mass or energy). The rate of change of the $L_2$ norm monitors transfer of energy between the subset of resolved scales and the rest of the resolved scales. When the rate of change of the $L_2$ norm exceeded a prescribed tolerance, mesh refinement was performed. The implementation of this criterion to perform mesh refinement yielded very encouraging results for the problem of detecting singularities for the inviscid Burgers and focusing Schr\"odinger equations. What we offer in the current work is a quantitative explanation of the reason behind the success of this approach to mesh refinement. In particular, we show that monitoring the rate of change of the $L_2$ norm allows us to monitor the rate of change of the error of the reduced model. Thus, by controlling the $L_2$ norm is a good way to  control the error.

In order for the algorithm to be successful one assumes the existence of an accurate reduced model. This guarantees that the transfer of activity to the unresolved scales is captured correctly so that there is not excessive or, worse, insufficient refinement. For the examples presented here we use a reduced model, called the $t$-model, which was derived through the Mori-Zwanzig formalism \cite{CHK3,HS06,S10-2}. In Sections \ref{burgers} and \ref{schrodinger} we provide convergence proofs of the $t$-model for Burgers and critical focusing Schr\"odinger respectively, as long as the solutions of these equations remain smooth. These proofs will be needed to exhibit the connection between model reduction and mesh refinement. 

We should note that closer examination of the mathematical operations in the two proofs shows that there are some common features. This is due to the fact that the reduced model has the same structure. Also, both cases examined here are equations that conserve the $L_2$ norm of the solution of the full system, while the latter is still smooth. The common structure of the reduced model for the two cases suggests that the same proof would go through for other partial differential equations that share the same features with the ones investigated here. In particular, the proof goes through as is for the subcritical and supercritical focusing Schr\"odinger equations.

%%%%%%%%%%%%%%%%End of Introduction%%%%%%%%%%%%%%%%%%%%%

\section{The mesh refinement algorithm}\label{algorithm}

A successful mesh refinement algorithm detects accurately the cascade of activity through the available simulated scales. A successful reduced model reproduces accurately the cascade of activity from the resolved scales to the unresolved ones. The mesh refinement algorithm presented in \cite{S07} starts by splitting the available simulated physical scales into two parts, call them $F$ and $G$ (the notation will become clearer in the following sections). The scales in $F$ are larger than the scales in $G$ (smaller if one thinks in terms of wavenumbers). Then, one constructs a reduced model for the scales in $F$ which, hopefully, captures correctly the cascade of activity to the scales in $G$. Note that one cannot wait until there is transfer of activity out of the scales in $G$ to even smaller scales. One must act earlier than that to guarantee that the simulation remains well resolved. The proposed mesh refinement algorithm monitors the rate of change of the $L_2$ norm of the scales in $F.$ When a prescribed threshold is crossed, the algorithm adds to the simulation scales that are smaller than those in $G$. This process can be repeated until one reaches the limit of scales that can be simulated based on the available computational power or until the algorithm does not detect the need for further refinement.

While the algorithm has performed very well as shown in \cite{S07}, there is a need to explain why monitoring the rate of change of the $L_2$ norm is a reliable indicator of the error. We will do that by showing that, indeed, the rate of change of the $L_2$ norm accounts for the rate of change of the error of the reduced model. If there is no transfer of activity from $F$ to $G,$ then the rate of change of the $L_2$ norm and the rate of change of the error of the reduced model will be zero (down to the numerical precision used). So, monitoring the rate of change of the $L_2$ norm is indeed a good way to monitor the error and thus a reliable indicator of the need for mesh refinement.

%%%%%%%%%%%%%%%%End of Section 1%%%%%%%%%%%%%%%%%%%%%

\section{Two examples}\label{examples}

In this section we present convergence proofs for the $t$-model for the 1D Burgers and critical focusing Schr\"odinger equations. In both cases we assume periodic boundary conditions which allows us to expand the solution in Fourier series. 

\subsection{Burgers equation}\label{burgers}

The 1D Burgers equation is given by 
\begin{equation}\label{burgersequation}
u_t+u u_x = 0.
\end{equation}
Equation (\ref{burgersequation}) should be supplemented with an initial condition $u(x,0)=u_0(x)$ and boundary conditions. We solve (\ref{burgersequation}) in the interval $[0,2\pi]$ with periodic boundary conditions. This allows us to expand the solution in Fourier series
$$u_{M}(x,t )=\underset{k \in F \cup G}{\sum} u_k(t) e^{ikx},$$
where $F \cup G=[-\frac{M}{2},\frac{M}{2}-1].$ We have written the set of Fourier modes as the union of two sets 
in anticipation of the construction of the reduced model comprising only of the modes in $F=[-\frac{N}{2},\frac{N}{2}-1],$ where $ N < M.$
The equation of motion for the Fourier mode $u_k$ becomes
\begin{equation}
\label{burgersode}
 \frac{d u_k}{dt}=- \frac{ik}{2} \underset{p, q \in F \cup G}{\underset{p+q=k  }{ \sum}} u_{p} u_{q}.
\end{equation}
The $t$-model \cite{HS06} is given by 
\begin{gather}\label{burgerstmodelreno}
 \frac{d v_k}{dt}=- \frac{ik}{2} \underset{p, q \in F }{\underset{p+q=k  }{ \sum}} v_{p} v_{q} +  t \biggl( (-\frac{ik}{2}) \underset{q \in G}{\underset{p \in F}{\underset{p+q=k  }{ \sum}}}  v_{p}  \bigl(- \frac{iq}{2} \underset{r,s \in F }{\underset{r+s=q  }{ \sum}} v_{r} v_{s} \bigr)   + \\
  (-\frac{ik}{2}) \underset{q \in F}{\underset{p \in G}{\underset{p+q=k  }{ \sum}}}   \bigl(- \frac{ip}{2} \underset{r,s \in F }{\underset{r+s=qp }{ \sum}} v_{r} v_{s} \bigr) v_{q}  \biggr)  \notag
\end{gather} 
for $k \in F.$  

To proceed we will need the partial differential equivalent of \eqref{burgerstmodelreno}. Recall that we have set $N=\frac{M}{2}$ and define
$$P u_M(x,t)=u_{N}(x,t)=\underset{k \in F}{\sum} u_k(t) e^{ikx},$$ i.e. the projection of the solution on the $N$ resolved modes. Also, $$Qu_M(x,t)=(I-P) u_M(x,t)=\underset{k \in G}{\sum} u_k(t) e^{ikx}.$$ Let $u,U$ be two finite Fourier series and set $$B[u,U]=\p_x \frac{1}{2}(uU).$$
Let $v_k(t)$ satisfy \eqref{burgerstmodelreno} and define $$v(x,t)= \underset{k \in F}{\sum} v_k(t) e^{ikx}$$ and $$\Gamma(x,t)=-(I-P)B(v,v).$$ Then, equation \eqref{burgerstmodelreno} is equivalent to the PDE 
\begin{equation}\label{burgerstmodelpde}
\p_tv + PB(v,v) = -  tP [B[v,\Gamma]+B[\Gamma,v]]
\end{equation}
The energy decay for the $t$-model is given by \cite{HS06}
\begin{equation} \label{burgersenergydecay}
\frac{d}{dt} \frac{1}{2} ||v||_2^2 =-  t ||(I-P) B[v,v] ||_2^2=-  t ||\Gamma ||_2^2 ,
\end{equation}
where $$||f(x,t)||_2^2\equiv \int_0^{2\pi} f^2(x,t)dx.$$ We also define the inner product $(\cdot,\cdot)$ as
$$(f(x,t),g(x,t))\equiv  \int_0^{2\pi} f(x,t)g(x,t)dx.$$ 
Note that for the solution of the full system (involving modes in $F \cup G$) we have
\begin{equation} \label{burgersenergyconserve}
\frac{d}{dt} \frac{1}{2} ||u||_2^2 =0.
\end{equation}
We want to estimate rate of change of $||v-Pu||_2^2.$ We are interested in the time interval when the solution is still smooth so that all the mathematical operations we perform are valid. 

We begin by writing the projected form of \eqref{burgersequation},
\begin{equation}\label{projectedburgersequation}
\p_t Pu+PB[u,u] = 0.
\end{equation}
We subtract \eqref{burgerstmodelpde} and \eqref{projectedburgersequation}, multiply by $v-Pu$ and integrate over $[0,2\pi]$ to obtain
\begin{equation}\label{burgerssubtraction}
\frac{d}{dt}\frac{1}{2}||v-Pu||_2^2=-(v-Pu,PB(v,v-PB[u,u])-t(v-Pu,P[ B[v,\Gamma] + B[\Gamma,v] ]).
\end{equation}
By standard arguments we find 
\begin{equation}\label{burgerssubtractionmarkov}
-(v-Pu,PB[v,v]-PB[u,u]) \leq C_1 ||v-Pu||_2^2 +C_2 ||Qu||_2^2.
\end{equation}
For the term $-t(v-Pu,P[ B[v,\Gamma] + B[\Gamma,v] ])=-t(v-Pu,\p_x(v\Gamma))$ we find
\begin{gather}
-t(v-Pu,\p_x(v\Gamma))=-t(v,\p_x(v\Gamma))+t(Pu,\p_x(v\Gamma)) \notag \\
=t(\p_x v,v\Gamma)-t(\p_x Pu,v\Gamma)=t(\p_x \frac{v^2}{2},\Gamma)-t(\p_x Pu (v-Pu),\Gamma) \notag \\
-t(\p_x \frac{(Pu)^2}{2},\Gamma) \notag \\
\leq -t ||(I-P)\p_x \frac{v^2}{2}  ||_2^2 +t ||\p_x Pu  ||_{\infty} ||v-Pu ||_2 ||\Gamma ||_2 \notag \\
+t  ||(I-P)\p_x \frac{(Pu)^2}{2} ||_2 ||\Gamma ||_2  \notag \\
\leq -t ||\Gamma  ||_2^2 + \frac{t}{2} ||\p_x Pu  ||^2_{\infty} ||v-Pu ||^2_2 + \frac{t}{2}  ||\Gamma ||^2_2 \notag \\
+\frac{t}{2}  ||(I-P)\p_x \frac{(Pu)^2}{2} ||^2_2 + \frac{t}{2}||\Gamma ||^2_2.  \notag
\end{gather}
After cancellations, we get 
\begin{equation}\label{burgerssubtractionmemory}
-t(v-Pu,\p_x(v\Gamma)) \leq  \frac{t}{2} ||\p_x Pu  ||^2_{\infty} ||v-Pu ||^2_2  + \frac{t}{2}  ||(I-P)\p_x \frac{(Pu)^2}{2} ||^2_2.
\end{equation}
We combine \eqref{burgerssubtraction}, \eqref{burgerssubtractionmarkov} and \eqref{burgerssubtractionmemory} to find
\begin{gather}
\frac{d}{dt}\frac{1}{2}||v-Pu||_2^2 \leq (C_1+t C_3) ||v-Pu||_2^2 +C_2 ||Qu||_2^2 \notag \\
+ \frac{t}{2}  ||(I-P)B[Pu,Pu] ||^2_2 \label{burgerserror}
\end{gather}
where we have used $B[Pu,Pu]=\p_x \frac{(Pu)^2}{2}.$ The convergence proof can be completed using Gronwall's inequality.

By examining \eqref{burgersenergydecay} and the last term in \eqref{burgerserror}, we see that the rate of change of the $L_2$ norm for the modes in $F$ as given by the $t$-model has the same functional form as the contribution of the $t$-model term to the inhomogeneous part of the rate of error of the $t$-model. The difference is that in \eqref{burgersenergydecay} the term is evaluated using the $t$-model solution while in \eqref{burgerserror} it is evaluated using the projection of the solution of the full system. Thus, the rate of change of the $L_2$ norm for the reduced model is related to the rate of change of the error for this model. 

In fact, this provides us with two alternative ways to monitor the cascade of activity from $F$ to $G.$ The first is to solve the full system and the reduced model simultaneously and monitor $ t ||(I-P) B(v,v) ||_2^2$ (see \eqref{burgersenergydecay}). When this exceeds a prescribed tolerance, we can add (if available) more modes and continue. The second way is to solve the full system and use the modes in $F$ to compute $ t ||(I-P)B[Pu,Pu] ||^2_2.$ When this exceeds a prescribed tolerance we can add (if available) more modes and continue. 

%%%%%%%%%%%%%%%%End of Section 2%%%%%%%%%%%%%%%%%%%%%

\subsection{Critical focusing Schr\"odinger equation}\label{schrodinger}

The 1D critical focusing Schr\"odinger equation is given by 
\begin{equation}\label{schrodingerequation}
iu_t+\Delta u+|u|^4u = 0.
\end{equation}
Equation (\ref{schrodingerequation}) should be supplemented with an initial condition $u(x,0)=u_0(x)$ and boundary conditions. We solve (\ref{schrodingerequation}) in the interval $[0,2\pi]$ with periodic boundary conditions. This allows us to expand the solution in Fourier series
$$u_{M}(x,t )=\underset{k \in F \cup G}{\sum} u_k(t) e^{ikx},$$
where $F \cup G=[-\frac{M}{2},\frac{M}{2}-1].$ We have written the set of Fourier modes as the union of two sets 
in anticipation of the construction of the reduced model comprising only of the modes in $F=[-\frac{N}{2},\frac{N}{2}-1],$ where $ N < M.$
The equation of motion for the Fourier mode $u_k$ becomes
\begin{equation}
\label{schrodingerode}
 \frac{d u_k}{dt}=-i k^2 u_k+ i \underset{k_1 \ldots, k_5 \in F \cup G}{\underset{k_1-k_2+k_3-k_4+k_5=k  }{ \sum}} u_{k_1} u^{*}_{k_2} u_{k_3} u^{*}_{k_4} u_{k_5}
\end{equation}
where $u^*_k$ is the complex conjugate of $u_k.$

The $t$-model \cite{S10-2} is given by (for $k \in F$)
\begin{multline}\label{schrodingertmodelreno}
\frac{d}{dt}v_k=-i k^2 v_k+ i \underset{k_1 \ldots, k_5 \in F}{\underset{k_1-k_2+k_3-k_4+k_5=k  }{ \sum}} v_{k_1} v^{*}_{k_2} v_{k_3} v^{*}_{k_4} v_{k_5} \\
+t \biggr[ 3i\underset{k_1 \in G  ,\, k_2,\ldots,k_5 \in F}{\underset{k_1-k_2+k_3-k_4+k_5=k  }{ \sum}}    R_{k_1}(\hat{v}) v^{*}_{k_2} v_{k_3} v^{*}_{k_4} v_{k_5}  \\
+2i \underset{k_2 \in G ,\, k_1,k_3,\ldots,k_5 \in F }{\underset{k_1-k_2+k_3-k_4+k_5=k  }{ \sum}}  v_{k_1}R^{*}_{k_2}(\hat{v})  v_{k_3} v^{*}_{k_4} v_{k_5} \biggl] 
\end{multline}
where $\hat{v}=(v_k)_{k \in F}$ and $R_k(v)=  i \underset{k_1 \ldots, k_5 \in F }{\underset{k_1-k_2+k_3-k_4+k_5=k  }{ \sum}} v_{k_1} v^{*}_{k_2} v_{k_3} v^{*}_{k_4} v_{k_5}.$

To proceed we will need the partial differential equivalent of \eqref{schrodingertmodelreno}. Recall that we have set $N=\frac{M}{2}$ and define
$$P u_M(x,t)=u_{N}(x,t)=\underset{k \in F}{\sum} u_k(t) e^{ikx},$$ i.e. the projection of the solution on the $N$ resolved modes. Also, $$Qu_M(x,t)=(I-P) u_M(x,t)=\underset{k \in G}{\sum} u_k(t) e^{ikx}.$$ We set $$B[z_1,z_2,z_3,z_4,z_5]=z_1 z^*_2 z_3 z^*_4 z_5.$$
Let $v_k(t)$ satisfy \eqref{schrodingertmodelreno} and define $$v(x,t)= \underset{k \in F}{\sum} v_k(t) e^{ikx}$$ and $$\Gamma(x,t)=i(I-P)B(v,v,v,v,v).$$ Then, equation \eqref{schrodingertmodelreno} is equivalent to the PDE 
\begin{equation}\label{schrodingertmodelpde}
\p_tv = i \Delta v +i PB[v,v,v,v,v] + i3t P B[\Gamma,v,v,v,v] +i2t PB[v,\Gamma,v,v,v].
\end{equation}
The energy decay for the $t$-model is given by \cite{S10-2}
\begin{equation} \label{schrodingerenergydecay}
\frac{d}{dt} \frac{1}{2} ||v||_2^2 =-  t ||(I-P) B[v,v,v,v,v] ||_2^2=-  t ||\Gamma ||_2^2,
\end{equation}
where $$||f(x,t)||_2^2\equiv \int_0^{2\pi} |f(x,t)|^2dx.$$ We also define the inner product $(\cdot,\cdot)$ as
$$(f(x,t),g(x,t))\equiv  \int_0^{2\pi} f^*(x,t)g(x,t)dx.$$
Note that for the solution of the full system (involving modes in $F \cup G$) we have
\begin{equation} \label{schrodingerenergyconserve}
\frac{d}{dt} \frac{1}{2} ||u||_2^2 =0.
\end{equation}
We want to estimate rate of change of $||v-Pu||_2^2.$ We are interested in the time interval when the solution is still smooth so that all the mathematical operations we perform are valid. 

We begin by writing the projected form of \eqref{schrodingerequation},
\begin{equation}\label{projectedschrodingerequation}
i\p_t Pu+\Delta Pu+PB[u,u,u,u,u] = 0.
\end{equation}
We subtract \eqref{schrodingertmodelpde} and \eqref{projectedschrodingerequation}, multiply by $(v-Pu)^*$ and integrate over $[0,2\pi]$ to obtain
\begin{gather}
\frac{d}{dt}\frac{1}{2}||v-Pu||_2^2=  Re \bigg[ (v-Pu, \p_t (v-Pu) )\bigg]= \notag \\
Re \bigg[ (v-Pu,i\Delta (v-Pu)) \bigg ] \notag \\ 
+Re \bigg [ i(v-Pu,PB[v,v,v,v,v]-PB[u,u,u,u,u]) \bigg ]\notag \\
+Re \bigg[ (v-Pu,i3tPB[\Gamma,v,v,v,v]) \bigg] \notag \\
+Re \bigg[ (v-Pu,i2tPB[v,\Gamma,v,v,v]) \bigg]. \label{schrodingersubtraction}
\end{gather}
For the term $Re \bigg[ (v-Pu,i\Delta (v-Pu)) \bigg ],$ integration by parts yields
$$ Re \bigg[ (v-Pu,i\Delta (v-Pu)) \bigg ] =0.$$
Also, by standard arguments, we find
\begin{equation}\label{schrodingersubtractionmarkov}
Re \bigg [ i(v-Pu,PB[v,v,v,v,v]-PB[u,u,u,u,u]) \bigg ] \leq C_1 ||v-Pu||_2^2 +C_2 ||Qu||_2^2 .
\end{equation}
We proceed with the terms $Re \bigg[ (v-Pu,i3tPB[\Gamma,v,v,v,v]) \bigg]$ and $Re \bigg[ (v-Pu,i2tPB[v,\Gamma,v,v,v]) \bigg].$ We have 
\begin{gather}
Re \bigg[ (v-Pu,i3tPB[\Gamma,v,v,v,v]) \bigg]=\underset{I}{\underbrace{Re \bigg[ i3t(v,PB[\Gamma,v,v,v,v]) \bigg] }}  \notag \\
+\underset{II}{\underbrace{Re \bigg[ -i3t(Pu,PB[\Gamma,v,v,v,v]) \bigg] }}  \label{one}
\end{gather}
and
\begin{gather}
Re \bigg[ (v-Pu,i2tPB[v,\Gamma,v,v,v]) \bigg]=\underset{III}{\underbrace{Re \bigg[ i2t(v,PB[v,\Gamma,v,v,v]) \bigg] }} \notag \\
+\underset{IV}{\underbrace{Re \bigg[ -i2t(Pu,PB[v,\Gamma,v,v,v]) \bigg] }}.  \label{two}
\end{gather}
From \eqref{schrodingerenergydecay}, we find 
\begin{equation}\label{schrodingersubtractionmemory}
I+III=-  t ||\Gamma ||_2^2.
\end{equation}
To obtain estimates for $II$ and $IV$ is more involved but it follows the same procedure as the one used to obtain the estimate \eqref{burgerssubtractionmemory} for Burgers equation. The fact that for Schr\"odinger we have a nonlinear term of higher degree leads to more terms but other than that the procedure is straightforward. 

We have for the term $II$ (omitting the $Re$ part and the multiplicative factor) 
\begin{gather}
(Pu,PB[\Gamma,v,v,v,v])=(Pu,PB[\Gamma,v-Pu,v-Pu,v-Pu,v-Pu]) \notag \\
+2(Pu,PB[\Gamma,v-Pu,v-Pu,v-Pu,Pu]) \notag \\
+2(Pu,PB[\Gamma,Pu,v-Pu,v-Pu,v-Pu]) \notag \\
+4(Pu,PB[\Gamma,Pu,Pu,v-Pu,v-Pu]) \notag \\
+(Pu,PB[\Gamma,v-Pu,Pu,v-Pu,Pu]) \notag \\
+(Pu,PB[\Gamma,Pu,v-Pu,Pu,v-Pu]) \notag \\
+2(Pu,PB[\Gamma,Pu,v-Pu,Pu,Pu]) \notag \\
+2(Pu,PB[\Gamma,Pu,Pu,v-Pu,Pu]) \notag \\
+(Pu,PB[\Gamma,Pu,Pu,Pu,Pu]) \notag
\end{gather}
By the Cauchy-Schwarz inequality, the Sobolev inequality on the torus \cite{benyi} and the assumption of smoothness of the solution we get
\begin{gather}
Re \bigg[ -i3t(Pu,PB[\Gamma,v,v,v,v]) \bigg] \leq tC_3 (N)[||v-Pu||_2^2]^4 \notag \\
+tC_4(N)  [||v-Pu||_2^2]^3  + tC_5(N) [||v-Pu||_2^2]^2+tC_6 ||v-Pu||_2^2 \notag \\
+1152t||(I-P)B[Pu,Pu,Pu,Pu,Pu]||_2^2 +\frac{1}{2} t ||\Gamma||_2^2 \label{three}
\end{gather}
Similarly,
\begin{gather}
Re \bigg[ -i2t(Pu,PB[\Gamma,v,v,v,v]) \bigg] \leq tC_7 (N)[||v-Pu||_2^2]^4 \notag \\ 
+tC_8(N)  [||v-Pu||_2^2]^3  + tC_9(N) [||v-Pu||_2^2]^2+tC_{10} ||v-Pu||_2^2 \notag \\
+512t||(I-P)B[Pu,Pu,Pu,Pu,Pu]||_2^2 +\frac{1}{2} t ||\Gamma||_2^2 \label{four}
\end{gather}
We combine \eqref{schrodingersubtraction}-\eqref{four} and we find
\begin{gather}
\frac{d}{dt}\frac{1}{2}||v-Pu||_2^2 \leq tC_1 (N)[||v-Pu||_2^2]^4 +tC_2(N)  [||v-Pu||_2^2]^3  \notag \\
+ tC_3(N) [||v-Pu||_2^2]^2+(C_4+tC_5) ||v-Pu||_2^2 + C_6 ||Qu||_2^2 \notag \\
+1664t||(I-P)B[Pu,Pu,Pu,Pu,Pu]||_2^2  \label{schrodingererror}
\end{gather}
for some constants $C_1(N),C_2(N),C_3(N),C_4,C_5$ and $C_6.$ Note that the use of the Sobolev inequality and a crude upper bound estimate leads to the dependence of the constants $C_1(N),C_2(N),C_3(N)$ on the number $N$ of modes in $F.$ In particular, as $N \rightarrow \infty$ these constants grow polynomially in $N.$ This does not present a problem for the convergence proof because, after integration of \eqref{schrodingererror},  these constants end up multiplying the inhomogeneous terms. As long as the solution is smooth (as we have assumed here), these inhomogeneous terms vanish fast enough so that convergence is achieved for all time instants. 

The numerical factor 1664 in front of $t||(I-P)B[Pu,Pu,Pu,Pu,Pu]||_2^2$ arises because we have arranged the contributions of the different terms to $t ||\Gamma||_2^2$ so that they cancel the $-t ||\Gamma||_2^2$ term (see \eqref{schrodingersubtractionmemory}). This is done so that the in the final estimate all the nonhomogenous terms involve only expressions depending on the solution of the full system. This allows the convergence proof to be completed by using inequalities for differential equations (see e.g. Section 2.4 in \cite{pachpatte}). 

As in the case of Burgers equation, by examining \eqref{schrodingerenergydecay} and the last term in \eqref{schrodingererror}, we see that the rate of change of the $L_2$ norm for the modes in $F$ as given by the $t$-model has the same functional form as the contribution of the $t$-model term to the inhomogeneous part of the rate of error of the $t$-model. In \eqref{schrodingerenergydecay} the term is evaluated using the $t$-model solution while in \eqref{schrodingererror} it is evaluated using the projection of the solution of the full system. Thus, the rate of change of the $L_2$ norm for the reduced model is related to the rate of change of the error for this model.

%%%%%%%%%%%%%%%%End of Section 3%%%%%%%%%%%%%%%%%%%%%

\section{Discussion and future work}
We have presented a novel approach for mesh refinement which is based on the use of a reduced model. For specific cases, we have shown the mathematical reason which allows us to use a reduced model to detect accurately the need for mesh refinement. 

In the current work the reduced model was used to detect the cascade of activity across Fourier modes (from large to small scales). The same idea can be applied in a non-periodic setting by using alternative spectral expansions or non-spectral methods. Note that for the case of spectral expansions we offer only a way to add modes in Fourier space. This leads to uniform refinement in real space. On the other hand, non-spectral methods, like finite difference or finite volume methods, allow to detect also {\it where} in space to define in addition to {\it when}. One way to do that is to divide the domain in subdomains and apply the mesh refinement algorithm on the subdomains. As in the current work, to be able to apply the algorithm requires that one has a good reduced model before using it to detect the need for mesh refinement. This is ongoing work and will be presented in a future publication.

%%%%%%%%%%%%%%%%End of Discussion%%%%%%%%%%%%%%%%%%%%%

\section*{Acknowledgements} 
I would like to thank Dr. J. Li for very helpful discussions and comments. This work was supported by the Applied Mathematics Program within the Department of Energy (DOE) Office of Advanced Scientific Computing Research (ASCR) as part of the Collaboratory on Mathematics for Mesoscopic Modeling of Materials (CM4).


\begin{thebibliography}{99}


\bibitem{almgren}
Almgren A.S., Bell J.B., Colella P., Howell L.H. and Welcome M.L., A conservative adaptive projection method for the variable-density incompressible Navier-Stokes equations, J. Comp. Phys. 142 (1998) pp. 1-46.

\bibitem{benyi}
Benyi A. and Oh T., The Sobolev inequality on the torus revisited, Publ. Math. Debrecen 83 (2013), no. 3, pp. 359--374.


\bibitem{berger1}
Berger, M. and Kohn, R., A rescaling algorithm for the numerical calculation of blowing-up solutions, 
Comm. Pure Appl. Math. 41 (1988) pp. 841-863.


\bibitem{berger2}
Berger M. and Colella P., Local adaptive mesh refinement for shock hydrodynamics, J. Comp. Phys. 82 (1989) pp. 62-84.
%
%\bibitem{bernstein}
%Bernstein D., Optimal prediction of Burger's equation, Multi. Mod. Sim. 6 (2007) pp. 27-52.
%
%
%\bibitem{binney}
%Binney J., Dowrick N., Fisher A., Newman M., The Theory of Critical Phenomena (An Introduction to 
%the Renormalization Group), The Clarendon Press, Oxford, 1992.
%
%\bibitem{boyd}
%Boyd J.P., Chebyshev and Fourier Spectral Methods, Dover, New York, 2001.
%
\bibitem{budd}
Budd, C. J., Huang, W. and Russell, R. D., Moving mesh methods for problems with
blow-up. SIAM Jour. Sci. Comput. 17 (1996) pp. 305-327.

\bibitem{ceniceros}
Ceniceros H.D. and Hou T.Y., An efficient dynamically adaptive mesh for potentially singular solutions, J. Comp. Phys. 172 (2001) pp. 609-639.

%
\bibitem{CHK00}
Chorin, A.J., Hald, O.H. and Kupferman, R.,
Optimal prediction and the Mori-Zwanzig representation of irreversible
processes. Proc. Nat. Acad. Sci. USA 97 (2000) pp. 2968-2973.
%
\bibitem{CHK3}
Chorin, A.J., Hald, O.H. and Kupferman, R., Optimal prediction with memory, Physica D 166 (2002) pp. 239-257.
%
\bibitem{CS05}
Chorin, A.J. and Stinis, P., Problem reduction, renormalization and memory,
Comm. App. Math. Comp. Sci. 1 (2005) pp. 1-27.
%
\bibitem{givon}
Givon, D., Kupferman, R. and Stuart, A., Extracting macroscopic 
dynamics: model problems and algorithms, Nonlinearity 17 (2004) pp. R55-R127.
%
%\bibitem{goldenfeld}
%Goldenfeld, N., Lectures on Phase Transitions and the Renormalization Group, Perseus Books, Reading, Mass., 1992.
%
%\bibitem{hair}
%Hairer, E.,  N\"orsett, S.E.  and Wanner, G., Solving Ordinary Differential Equations I, Springer, NY, 1987.
%
\bibitem{HS06}
Hald O.H. and Stinis P., Optimal prediction and the rate of decay for solutions of the Euler equations in two and three dimensions, Proc. Natl. Acad. Sci.,104, no. 16 (2007) pp. 6527-6532.
%
\bibitem{landman}
Landman, M.J.,  Papanicolaou, G.C., Sulem, C. and Sulem, P.,  Rate of blowup for solutions of the nonlinear Schršdinger equation at critical dimension, Phys. Rev. A 38 (1988) pp. 3837-3847.

\bibitem{pachpatte}
Pachpatte, B.G., Inequalities for Differential and Integral Equations, Academic Press, San Diego, 1998.

\bibitem{S07}
Stinis P., A phase transition approach to detecting singularities of PDEs, Comm. App. Math. Comp. Sci. Vol. 4 (2009), No. 1, 217-239.

\bibitem{S10-2}
Stinis P., Numerical computation of solutions of the critical nonlinear Schr\"odinger equation after the singularity, Multiscale Modeling and Simulation 10 (2012), pp. 48-60.

%\bibitem{sulem}
%Sulem C. and Sulem P.-L., The nonlinear Schr\"odinger equation - Self-focusing and wave collapse, Applied Mathematical Sciences, 139, Springer, New York, 1999.
%
%\bibitem{weinberg}
%Weinberg S., Why the renormalization group is a good thing, Asymptotic Realms of Physics: Essays in Honor of Francis E. Low, MIT Press, Cambridge MA, 1983.
%
%\bibitem{wilson}
%Wilson K.,  The renormalization group and critical phenomena, Rev. Mod. Phys. 55 (1983) pp. 583-600.
%
%\bibitem{zakharov}
%Zakharov V.E., Collapse of self-focusing Langmuir waves, Handbook of plasma physics 2, North Holland, Amsterdam, 1984.

\end{thebibliography}
\end{document}